\documentclass[oneside,11pt,a4paper]{article}

\usepackage{amsmath,amsfonts,amssymb,epsfig}

  \title{Hausdorff dimension of the set of extreme points of a
  self-similar set}
 \author{Andrew Tetenov and Ivan Davydkin , }
 \date{}

\begin{document}
\maketitle
\vspace{-1cm}
\begin{center}
\it (Gorno-Altaisk state university)
\end{center}

 \newcommand \tg {\mathop{\rm tg}\nolimits}
 \newcommand \ctg {\mathop{\rm ctg}\nolimits}
 \newcommand \arctg {\mathop{\rm arctg}\nolimits}
 \newcommand \arcctg {\mathop{\rm arcctg}\nolimits}
 \newcommand \re {\mathop{\rm Re}\nolimits}
 \newcommand \im {\mathop{\rm Im}\nolimits}
 \newcommand \rr {\mathbb{R}}
 \newcommand \diam {\mathop{\rm diam}\nolimits}

 \newtheorem{teo}{\sc Theorem}
 \newtheorem{sled}[teo]{\sc Corollary}
 \newtheorem{lem}[teo]{\sc Lemma}
 \newtheorem{prop}[teo]{\sc Proposition}
 \newtheorem{defin}[teo]{\sc Definition}
 \newcommand{\e}{\varepsilon}

 \begin{abstract}
If the system
 of contracting similitudes on $\mathbb R^2$ satisfies open convex set
 condition,
then the set of extreme points of the convex hull $\tilde{K}$ of
it's invariant self-similar set
 has Hausdorff dimension 0 .
   If, additionally,  all the  rotation angles
 $\alpha_i$ of the similitudes
 $\varphi_i$ are commensurable with $\pi$,
  then the set  $\tilde{K}$ is a convex polygon.
\end{abstract}

\section{Introduction}
 Let $\varphi_1,\ldots,\varphi_n$ -- be   contracting similitudes
 on
 $\rr^2$
 with  rotation angles $\alpha_i$.
 Let $K=\bigcup\limits^n_{i=1}\varphi_i(K)$ -- be  compact self-similar set
 with respect to  the system $\{\varphi_1,\ldots,\varphi_n\}$.

 Let $\tilde{K}$ --- be the convex hull   $H(K)$ of the set  $K$.

 The boundary $\partial\tilde{K}$ of the set $\tilde{K}$ may be
 viewed as consisting of two parts:
 a compact set  $F$ of all extreme points of
 $\tilde{K}$, which we call {\it vertices} of
 $\tilde{K}$ , and of an union of finite or
 countable family  of open line intervals
 $l_i$, called
  {\it the sides } of $\tilde{K}$.

 There is a number of well-known cases when the set $\tilde{K}$
 is a convex finite-sided polygon.
 In Example 1 we show how mappings $\varphi_i$ with irrational
 rotation angles $\alpha_i$ may produce self-similar sets
 $K$,  whose convex hull  $\tilde{K}$ has infinite set of sides
  (  and hence of vertices)  .
 So the question arises , what conditions must be imposed upon
  $\{\varphi_1,\ldots,\varphi_n\}$ to obtain the set $\tilde{K}$
  with finite number of sides and vertices and
   what is the structure of
   $\partial\tilde{K}$ in the case when
  the set $F$ or $\{l_i\}$ is infinite.

To give a partial answer to these questions, we formulate  open
 convex set condition (Definition \ref{osc}),
 implying some finiteness properties for
 $\partial\tilde{K}$.
Then we show  that, if the system
 ($\varphi_1,\ldots,\varphi_n$) satisfies open convex set
 condition,
then the set of vertices of
 $\tilde{K}$ has Hausdorff dimension 0 (Theorem  \ref{t_1}).
   If, additionally,  all the  rotation angles
 $\alpha_i$ of the maps
 $\varphi_i$ are commensurable with $\pi$,
  then the set of vertices of $\tilde{K}$
is finite, so the set $\tilde{K}$ is a convex polygon
 (Corollary \ref{sz}).

 The authors express their appreciation to V.V.Aseev who suggested
 the problem.

 \section{Preliminaries.}

 Let $S$ be a system of contraction similitudes $\varphi_1,\ldots,\varphi_n$,
 (written in the complex form as  $\varphi_i(z)=q_i(z-z_i)e^{i\alpha_i}+z_i$, $0<q_i \leq q<1$,
 $\alpha_i\in[0,2\pi]$). The system $S$  defines
 the Hutchinson transformation
 $T(A)=\bigcup \varphi_i(A)$ on the space
 $C(\rr^2)$ of compact non-empty subsets of $\mathbf{R}^2$. The transformation
 $T$ is a contraction map on the space  $C(\rr^2)$ in the Hausdorff metrics.
  Compact invariant set $K=\bigcup\limits^n_{i=1}\varphi_i(K)$
   is a fixed point of the transformation $T$.

By  $H(X)$ we denote the convex hull of a set  $X$.

We shall denote by $\tilde{K}$ the convex hull $H(K)$ of the
invariant set $K$ with respect to the system $S$.

 Using notation
similar to that of  \cite{Hut}, we denote\\
 $\varphi_{i_1\ldots i_m}\!\!=\!\varphi_{i_1}\!\!\circ\varphi_{i_2}\!\!\circ
 \!\ldots\!\circ\varphi_{i_m}(K)$,\\
 $K_{i_1\ldots i_m}=\varphi_{i_1\ldots i_m}(K)$,
 $\tilde{K}_{i_1\ldots i_m}=\varphi_{i_1\ldots i_m}(\tilde{K})$.\\
By $Q_{i_1\ldots i_m}$ we denote $\tilde{K}_{i_1\ldots i_m}\cap
\partial\tilde{K}$   and\\  by
$F_{i_1\ldots i_m}\!=\!\tilde{K}_{i_1\ldots i_m}\cap F$-- the set
of extreme points belonging to $\tilde{K}_{i_1\ldots i_m}$.\\
 Let  $\tilde{C}(\rr^2)$  be the set of all \underline{convex} compact
 subsets in  $\rr^2$.  Define the transformation $\tilde{T}$ on the space
 $\tilde{C}(\rr^2)$ by:
 $\tilde{T}(A)=H(T(A))=H(\cup \varphi_i(A))$ for any
$A\subset \tilde{C}(\rr^2)$.

 \begin{prop}\label{prop_1}.
Let  $S=\{\varphi_1,\ldots,\varphi_n\}$ ---
 be a set of contracting similitudes on $\rr^2$
 with $Lip(f_i)=q_i$, $0< q_i\leq q<1$. The transformation $\tilde{T}$
 is a con\-trac\-tion map of the space $\tilde{C}(\rr^2)$ in Hausdorff metrics,
 and the set $\tilde{K}$ is it's fixed point.
 \end{prop}
 {\it Proof.} For any two compact sets $U_1$, $U_2$ the Hausdorff distance
 between their convex hulls $d(H(U_1),H(U_2))$ is less or equal then $
 d(U_1,U_2)$. Since  $\tilde{T}(A)=H(T(A))$, we obtain that
  $d(\tilde{T}(A),\tilde{T}(B))\leq
 d(T(A),T(B))<qd(A,B)$. 

 The equality $\tilde{T}(\tilde{K}) =H(\cup \varphi_i(\tilde{K}))=
H(\cup\varphi_i(K))=\tilde{K}$,  shows that  $\tilde{T}$ leaves
$\tilde{K}$
 fixed.
 \vspace{10pt}

 {\bf Example 1.}

 Consider a system
   $\varphi_1(z)=\frac{1}{3}ze^{i\alpha}-1$,
 $\varphi_2(z)=\frac{1}{3}ze^{i\alpha}+1$, where  $\alpha=r\pi $
 with irrarional r.
Let us show that convex invariant set $\tilde{K}$
 with respect to  the system $\varphi_1,\varphi_2$
  has infinite set of sides.
  As it follows from the Proposition \ref{prop_1},
  for any compact set $A\subset \rr^2$,
 $\lim\limits_{n\rightarrow\infty}\tilde T^n(A)=\tilde{K}$.

Let $P_0=\{\bar{0}\}$ and $P_n=\tilde{T}(P_0)$. Obviously, each
$P_n$ is a convex polygon and $P_0\subset P_1\subset \ldots
\subset P_n\subset\ldots$, so $\tilde{K}=\lim\limits_{n\rightarrow
+\infty}P_n=\bigcup P_n$.

 It's easy to see that $n$-th polygon
$P_n$ has $n$ pairs of opposite sides $l_k'$, $l_k''$, each having
length $2/3^{k-1}$ and forming an angle $(k-1)\alpha$ with
horisontal axis:\\

\begin{figure}[h]
\begin{center}
\epsfig{file=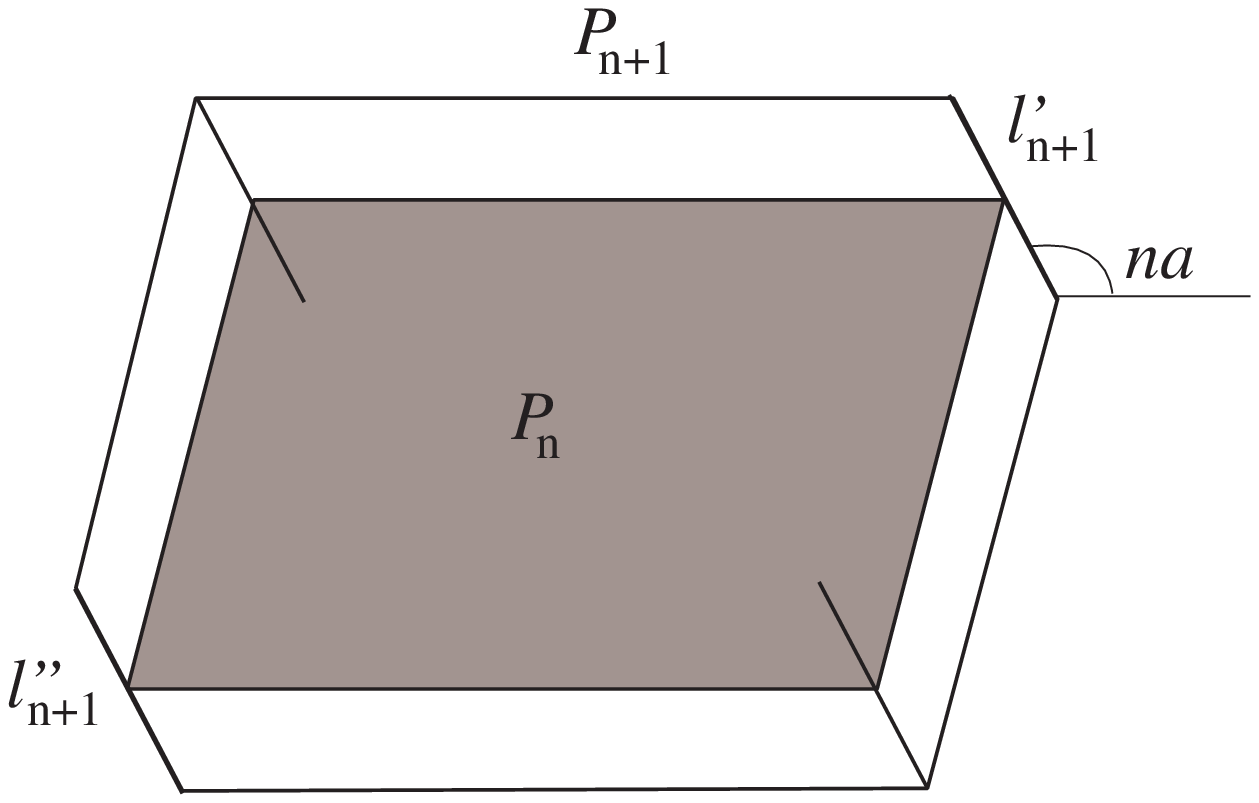,angle=0,scale=.5} \caption{} \label{f1}
\end{center}
\end{figure}
%

It's obvious for $n=2$. Assuming the fact is true for $P_n$,
observe that $P_{n+1}=P_n+l_{n+1}$, where $l_{k+1}$ is a line
segment $\left[-\frac{e^{ik\alpha}}{3^k},
\frac{e^{ik\alpha}}{3^n}\right]$ and $A+B$ denotes the set $\{a+b,
a\in A, b\in B\}$. Since neither of sides $l_k'$, $l_k''$ of $P_n$
is  parallel to $l_{n+1}$, $P_{n+1}$ has two more sides equal to
$l_{n+1}$. The sides $l_k'(P_n)$, $l_k''(P_n)$ are equal and
parallel for different $n$ and converge to the sides
$l_k'(\tilde{K})$ and $l_k''(\tilde{K})$ of $\tilde{K}$ having the
same length and direction.\\
 So $\tilde{K}$ has the sides $l_k'$
and $l_k''$ for each $k\in N$. Taking their sum we see that the
length of $\partial\tilde{K}$ is more or equal to
$4+\frac{4}{3}+\ldots+\frac{4}{3^n}+\ldots=6$.\\

\begin{figure}[h]
\begin{center}
\epsfig{file=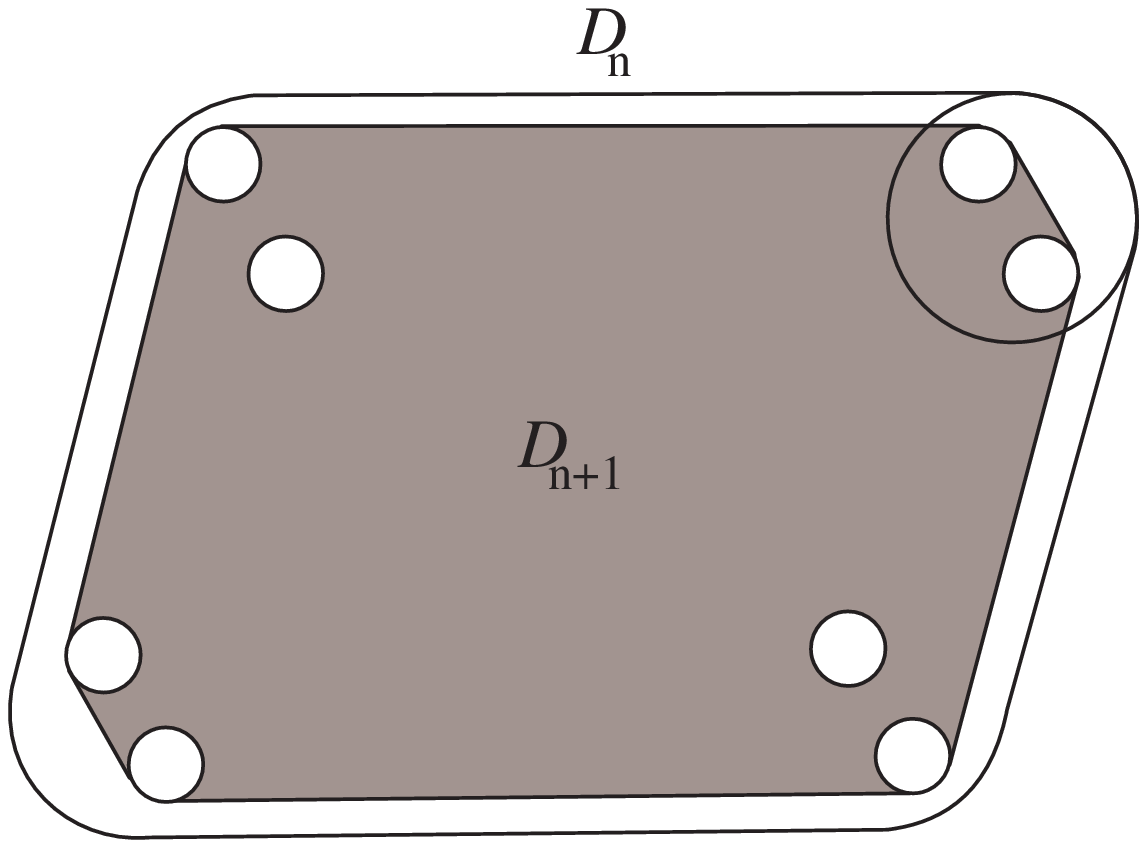,angle=0,scale=.5} \caption{} \label{f2}
\end{center}
\end{figure}
%
 To get the opposite inequality,
 consider the sets $D_n=\tilde{T}^n(D)$, where $D=\{x^2+y^2\leq
1\}$ is an unit disc of radias $\rho > 3$. They form a nested sequence
 $D_1\supset D_2\supset\ldots\supset D_n\supset\ldots$ and boundary $\partial
D_n$ of $D_n$ consists of $n$ pairs of sides  $l_k'$, $l_k''$,
each having length $2/3^{k-1}$ and forming an angle $(k-1)\alpha$
with horisontal axis and of $2n$ arcs of radius $\rho/3^{n+1}$, so it
has total length
$\frac{2\pi\rho}{3^{n+1}}+4+\frac{4}{3}+\ldots+\frac{4}{3^n}$. As n
tends to $\infty$ we obtain the opposite inequality
$H^1(\partial\tilde{K})\leq 6$, so $H^1(\partial\tilde{K})=6$. In
other words, 1-dimensional measure of $\partial\tilde{K}$ equals
the sum of  lengths of sides of $\tilde{K}$. Therefore the set of
vertices of $\tilde{K}$ has 1-dimensional Hausdorff measure 0.

\section{Open convex set condition.}

 \begin{defin}\label{osc}.
 A system $\{\varphi_1,\ldots,\varphi_n\}$
 satisfies open convex set condition (OCSC),
if there is  non-empty open convex set
 $O$, such that
 \begin{itemize}
 \item[(i)] $\cup\varphi_i(O)\subset O$,
 \item[(ii)]$\varphi_i(O)\cap\varphi_j(O)=\emptyset \; \mbox{\rm if } i\neq
 j$.
 \end{itemize}
 \end{defin}


 \begin{prop}\label{prop_2}.
 If the interior   $U$ of the set ${\tilde{K}}$ is non-empty, then
 \begin{itemize}
 \item[(i)] $\cup\varphi_i(U)\subset O$,
 \item[(ii)]$\varphi_i(U)\cap\varphi_j(U)=\emptyset$\  {\rm if}\  $i\neq j$.
 \end{itemize}
   \end{prop}
 {\it Proof.} The first inclusion is obvious. Let now $O$ be the open set
 from the OCS condition.
Since $\tilde{T}(\bar{O})\subset\bar{O}$, the set  $\tilde{K}$
 is contained in $\bar{O}$. Therefore
 $\varphi_i({\tilde{K}})\bigcap\varphi_j({\tilde{K}})\subset
 \varphi_i(\bar{O})\bigcap\varphi_j(\bar{O})\subset
 \varphi_i(\bar{O}\setminus O)  \bigcup\varphi_j(\bar{O}\setminus O)
 $. The latter set is closed and nowhere dense in  $\rr^2$ and the
 set
 $\varphi_i(U)\bigcap\varphi_j(U)$  is it's open subset. Therefore it is empty. 

 \begin{lem}\label{lem_1}.
 Let $U_1,U_2,...,U_n$ be closed convex domains
  with disjoint interiors and let $V$ be the convex hull of the
  set
 $\bigcup\limits_{i=1}^n U_k$. Let $Q_i=U_i\cap\partial V$,\
 $Q_{i,s}$ be the components of the set $Q_{i}$, and $m_i $
 be the number of such components.
 Then $\sum\limits_{i=1}^n m_i\leq 2n-2$.
 \end{lem}
{\it Proof.} The case $n=2$ is obvious, because $m_1=m_2=1$.

\begin{figure}[h]
\begin{center}
\epsfig{file=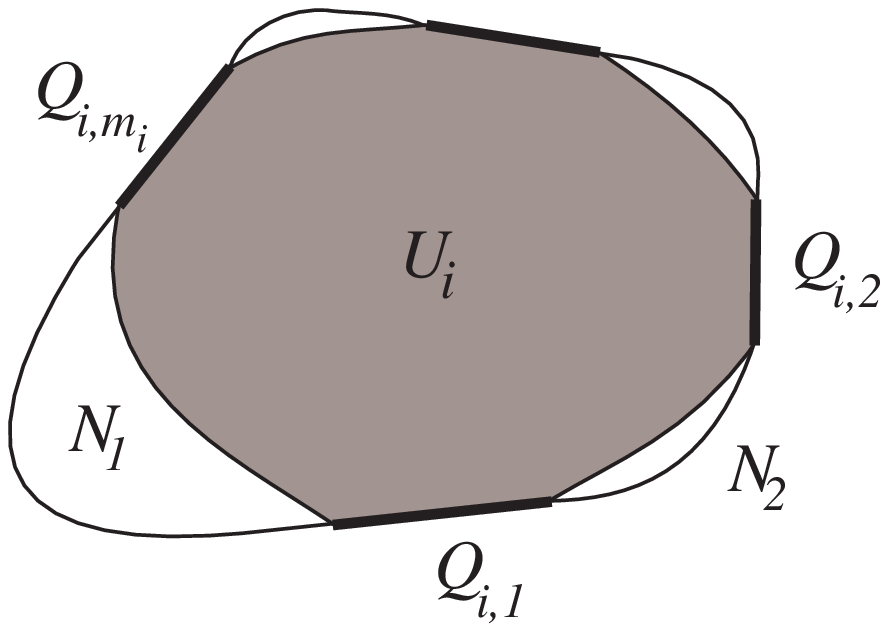,angle=0,scale=.5} \caption{} \label{f3}
\end{center}
\end{figure}

  Suppose it is already proved for any family consisting of
   $n-1\geq 2$  sets \ $U_i$. If
  all
 $m_i=1$, the statement also holds, so we can suppose that for
some $Q_i$, the number of it's components   $ m_i\geq 2$. Denote
the components of $V\setminus
 U_i $ by $N_k,\; k=1,\ldots, m_i$. Each of the sets  \ $U_j,\; j\neq i$
 is contained in the closure of some component\  $N_k$,
  and each $\bar{N_k}$
 contains at least one of the sets
  $U_j,\; j\neq i$: therefore $m_i\leq n-1$.

  Let $n_k $ be the number of all $U_j \subset \bar{N_k}$.
 Obviously ,\    $1\leq n_k \leq n-2 $.

Consider a family
  $\{U_i\}\cup \{U_j|U_j\subset\bar{N_k}\}$  and it's convex hull
 $\bar{N_k}\cup U_i$. The number of the sets $U_k$ in this family
 is
 $n_k+1 \leq n-1  $, so we can apply the statement of  Lemma
  and obtain
 $1+\sum\limits_{U_j\subset \bar{N_k}} m_j \leq 2n_k. $
 A total sum over all components $N_k$  gives the desired
 inequality
  $$ m= \sum\limits_{j\neq i} {m_j+m_i} \leq 2\sum n_k = 2(n-1).$$ 

 \begin{prop}\label{prop_3}.
 1) Each of the sets \space \space  $Q_i,i=1\ldots,n$
 contains finite number of components  $Q_{i,s}$,
  and their total number is less or equal to $2n-2$.\\
 2) If $n>2$ then for  $i\neq j$ the set $Q_{i,s}\bigcap Q_{j,t}$
 is either empty or  is a point .\\
 3) If a sequence of indices $(i_1,\ldots,i_n)$
 is an initial interval of $(j_1,\ldots,j_m)$, then
  $Q_{i_1,\ldots,i_n}\subset Q_{j_1,\ldots,j_m}$ .
 If neither of sequences $(i_1,\ldots,i_n)$, $(j_1,\ldots,j_m)$
 is an initial interval of the other
 then $Q_{i_1,\ldots,i_n}\cap Q_{j_1,\ldots,j_m}$
  is either empty or  is a point.\\
 4) If for some $j$,\space  $Q_j=\emptyset$,
 then for each sequence $(i_1,\ldots,i_n)$,\space $
 Q_{i_1,\ldots,i_n,j}=\emptyset$
 \end{prop}

 {\it Proof.}
The first statement results from the previous lemma.

To prove the second statement, suppose the set  $Q_{i,s}\cap
 Q_{j,t}$ contains two points $x,y $. Let $ l $ be a line segment
  with endpoints
 $x,y $. If $ l\subset Q_{i,s}\cap Q_{j,t} $ , then some
 half-neighbourhood of the point
 $\frac{x+y}{2}$ is contained in $ \tilde{K_i}\cap
 \tilde{K_j}$, which contradicts OCSC.
  By the very reason   $l\setminus\{x,y\}$
cannot be contained neither in $ \dot{\tilde{K_i}}$ nor in
$\dot{\tilde{K_j}}$. Since $n>2$,  $ l $ cannot be contained
   in the intersection of the boundaries of the sets
 $ \tilde{K_i}$ and $ \tilde{K_j}$ .

 The third statement follows directly from 2). 4) is obvious.

 \section{Finiteness of the set of sides of  $\tilde{K}$.}

 \begin{defin}\label{def_4}. We call a side $l$
  with endpoints $x_1,x_2\in F$ a side of  order 1,
  if there are  such $j_1\neq j_2$  that
 $x_1\in Q_{j_1},x_2\in Q_{j_2}$
 \end{defin}

 \begin{prop}\label{prop_5}.
 If  $\{\varphi_1,\ldots,\varphi_n\}$
 satisfies OCSC, then the set of sides of  order 1 is finite.
 \end{prop}
 {\it Proof.} The sides of order 1 are the closures of components
 of $\partial \tilde{K}\setminus  \bigcup\limits_{i,s} Q_{i,s}$ so from
 Proposition  \ref{prop_3} it follows that there is no more than
 $2n-2$ such sides.

 \begin{prop}\label{prop_6}.
 Each side  $l'\subset \partial \tilde{K}$  is an image
 $\varphi_{i_1i_2\ldots i_p}(l)$ of some side $l$ of  order 1.
 \end{prop}

 {\it Proof} Let $x',y'$ be the endpoints of the side  $l'$  .
  If $l'$ is contained in a component
 $Q_{i_1i_2\ldots i_p}$, then  $\varphi_{i_1i_2\ldots
 i_p}^{-1}(l')\subset \partial \tilde{K}$. It follows from the statement 3)
 of the Proposition \ref{prop_3} ,
  that among all components \ $Q_{i_1i_2\ldots i_{p'}}$,
containing the side
 $l'$, there is a component $Q_{i_1i_2\ldots i_p}$ of maximal order.
 There exist different
 $j,k =1,...,n$,  such that $x'\in Q_{i_1i_2\ldots i_pj}$,
 $y'\in Q_{i_1i_2\ldots
 i_pk}$. Then the side $l=\varphi_{i_1i_2\ldots i_p}^{-1}(l')$
is an side of order 1, because it's endpoints
 $x=\varphi_{i_1i_2\ldots i_p}^{-1}(x')$ and $y=\varphi_{i_1i_2\ldots
 i_p}^{-1}(y')$ belong to different components $Q_j$ and $Q_k$ .

 \begin{teo}\label{t_1}
If $\varphi_1,\ldots,\varphi_n$ satisfies OCSC and all the angles
$\alpha_1,\ldots,\alpha_n$ are commensurable with
 $\pi$, then the set of sides of $\tilde{K}$ is finite.
 \end{teo}

 {\it Proof}.  Let $\theta$ be the greatest common divisor
 of the angles   $\alpha_1,\ldots,\alpha_n,\pi$. Let $l_1,\ldots,l_m$
 be the edges of  order 1
, and $\beta_1,\ldots,\beta_m$ be the angles formed by
$l_1,\ldots,l_m$ and the  axis  $Ox$ . Each of the sides of
$\tilde{K}$ is an image $\varphi_{i_1\ldots i_p}(l_j)$ of some
side $l_j$ of order 1, and therefore the angle it forms with the
horizontal axis is of the form $\beta_j+n\theta$ . The set of all
such angles is finite , so the set of all sides of $\tilde{K}$ is
finite.  $\blacksquare$

 \subsection{Vertices and corner points.}

 Since every vertex $z_0$ of $\tilde{K}$ lies in
  $\partial\tilde{K}\cap K $  there is an infinite sequence
   $K_{i_1}\supset K_{i_1i_2}\supset\ldots
  K_{i_1i_2\ldots i_p}\supset\ldots$ such that
   $z_0= \bigcap\limits_{p=1}^\infty K_{i_1i_2...i_p}$.
So for each  $p$  there is a vertex $z_p=\varphi^{-1}_{i_1\ldots
i_p}(z_0)$. Thus each vertex $z_0$ has an infinite sequence of
predecessors $z_p \in\partial\tilde{K}$ such that
$z_0=\varphi_{i_1\ldots i_p}(z_p)$ for some $\varphi_{i_1\ldots
i_p}$.

   The vertex  is {\em periodic} if there is a periodic
   sequence $i_1,i_2,...i_p,...$ defining $z_0$. In this case
   one of $z_0$'s predecessors, say $z_m$, is a fixed point of
   some $\varphi_{i_1\ldots i_p}$.
   We call a vertex $z_0$ { \em a corner point} if right and left
   tangents to $\tilde{K}$ at $z_0$ do not coincide.

%

 \begin{prop}\label{prop_ugl}.
 Each corner point $z_0\in\partial\tilde K$  is a periodic vertex
 and both right and left tangents at a corner point $z_0$ are
 sides $l_0^+$and $l_0^-$.

 \end{prop}

 {\it Proof }. Let  $\theta_0$ be an angle between right and left
 tangents to $\tilde K$ at $z_0$, and
 $\theta_p$ an angle between right and left
 tangents to $\tilde K$ at $z_p$.
 Since  $\varphi_{i_1,..i_p}(z_p)=z_0$, the angle $\theta_p$ is more
 or equal to $\theta_0$ . The sum of $\theta_p$ over all of
 different predecessors of $z_0$ does not exceed $2\pi$ ,
 therefore the sequence $z_1,z_2,\ldots,z_p,\ldots$ contains no more
 than $2\pi/\theta_0$ different elements. This shows that $z_0$ is
 periodic.

Let now $V^+(z_0)$ denote a half-neighborhood  of the point $z_0$ in
$\partial\tilde K$ taken in positive direction. As we see from
Proposition  \ref{prop_3} for each $p$ there is unique $p$-tuple
$i_0,i_1,...,i_p$ for which $Q_{i_0i_1...i_p}\cup V^+(z_0)\setminus z_0$ is
non-empty for each $V^+(z_0)$. The sequence $Q_{i_0}\supset
Q_{i_0i_1}\supset...Q_{i_0i_1...i_p}\supset...$ defines  unique
sequence $z_1,z_2,\ldots,z_p,\ldots$ of predecessors of $z_0$
having the additional property that for each half-neighborhood
 $V^+(z_0)$ of $z_0$ and $V^+(z_p)$ of $z_p$ the intersection
$V^+(z_0)\setminus z_0\cap \varphi_{i_0i_1...i_p}(V^+(z_p))$ is non-empty.
Since the number of different $z_p$'s is finite, one of them, say $z_m$,
is a fixed point of some $ h=\varphi_{j_1\ldots j_q}$, satisfying
$h(V^+(z_p))\subset V^+(z_p)$ . The latter is possible
 only when $ V^+(z_p)$ is a straight line interval.


 \begin{prop}\label{prop_ugl}.
If all the angles  $\alpha_1,\ldots,\alpha_n$  are such that for
each set of non-negative integers  $k_1,\ldots,k_n,\;
k_1\alpha_1+\ldots
 +k_n\alpha_n\neq m\pi,m\in \mathbb{N}$,  then the set $F$
  of the vertices of $\tilde{K}$ is infinite.

 \end{prop}

 {\it Proof}. If the set of vertices of $\tilde{K}$ is finite, all
 of them are corner points . Since none of $\varphi_{i_1\ldots i_q}
 $ has rotation angle $m\pi$, $\partial\tilde{K}$ has no corner
 points, so the set $F$ is infinite.

  \subsection{The main theorem.}

 We call two systems   and
  $(\varphi_1',\ldots,\varphi_n')$
  {\it convex equivalent }, if they generate the same
  convex invariant set
   $\tilde{K}$.

   Let $S'= (\varphi_1',\ldots,\varphi_{n'}')$ be $p$-th refinement
   of the system   $S=(\varphi_1,\ldots,\varphi_n)$, i.e. the
   set of all mappings $\varphi_{i_1\ldots i_p}$ of order $p$.
   The system $S'$ is convex equivalent to $S$, and  $S'$ satisfies
   OCSC  if $S$ does.

 \begin{prop}.
For any system of similitudes $(\varphi_1,\ldots,\varphi_n)$
 satisfying the OCS condition, there is a
convex equivalent system
 $(\psi_1,\ldots,\psi_n)$
  satisfying the OCS condition
  such that
 \begin{equation} \label{e1}
  {\textstyle\rm\;all\; the\; sets\;}
   \psi_{i_1\ldots i_k}(\tilde{K}\cap\partial\tilde{K})
   {\textstyle\rm\;are\; connected\;}
 \end{equation}
 \begin{equation} \label{e2}
 \forall i \;\;\varphi_i(\tilde{K})\cap F=\emptyset
 \end{equation}
 \begin{equation} \label{e3}
 \forall i\neq j \;\;\varphi_i(\tilde{K})\cap F\not\subset\varphi_j(\tilde{K})
 \cap F.
 \end{equation}
 \end{prop}

 {\it Proof}. Suppose for some $i$ the set
  $Q_i=\tilde{K}_i\cap\partial\tilde{K}$, is not connected,
  and  $Q_{i,1} \ldots Q_{i,s}$
 are it's components. Let $\delta_i$ be  minimal distance
 between different components of  $Q_i$.
 Let $\delta$ be the smallest of all
 $\delta_i$ among non-connected $Q_i$'s.

 With each component $Q_{i,k}$ of $Q_i$ we associate a set
 $N_{i,k}$ of all unit outer normal vectors to $\partial\tilde K$
 at points $x\in Q_{i,k}$. If any edge $l\subset\partial\tilde K$
 has both of it's edges in $Q_i$ then it lies completely in one of
 it's components, say $Q_{i,k}$.

 Therefore for each $i$, the sets $N_{i,k}$ are disjoint closed
 arcs on the unit circle. Let $\theta_i$ be the  length of
 the shortest complementary arc to $\bigcup \limits_k N_{i,k}$.
 Let $\theta$ be the smallest  of all $\theta_i$'s.

Choose such  $p_0$ that $q^{p_0}\cdot\diam(\tilde{K})<\delta$.

 For any $p_0$-tuple $i_1\ldots i_{p_0}$ the set
 $\varphi_{i_1\ldots i_{p_0}}(
 \tilde{K})$ has diameter less than $\delta$, so the set
 $\varphi_{i_1\ldots i_{p_0}}(\tilde{K})$ may have  non-empty
 intersection with at most one of the components of $Q_{i_1}$.
 Therefore $Q_{i_1\ldots i_{p_0}}$ lies completely in some
 component of $Q_{i_1}$. Take any sequence $(j_1\ldots j_q)$.
 By  the same reason, for any $k\leq q$ the set
 $Q_{j_1\ldots j_p i_1\ldots i_{p_0}}$ must be contained
 completely in some component of $Q_{j_1\ldots j_k}$.
 Consider $S'= (\varphi_1',\ldots,\varphi_{n'}')$ be $p$-th refinement
   of the system   $S=(\varphi_1,\ldots,\varphi_n)$. By the above
   argument each non-empty set $Q_{ij}=\psi_{ij}(\tilde{K})\cup
   \partial\tilde{K}$ is contained in an unique component of the
   set $Q'_{ij} =\psi_i(\tilde{K})\cup
   \partial\tilde{K}$.
Replacing $S$ by $S'$  if necessary, we may suppose from that
moment
  that for each $i$, $j$, $Q_{ij}$ is
 contained completely in some component of $Q_i$.

 Assume some $Q_{i_1\ldots i_p}$ is non-connected
 and $Q^{(1)}_1$ and $Q^{(1)}_2$
 are   it's two adjoining components, joined by an arc
 $\Delta^{(1)}$ in $\partial \tilde{K}$
 with endpoints $\xi^{(1)}\in Q^{(1)}_1$
 and $\eta^{(1)}\in Q^{(1)}_2$. Then for each $1\leq k \leq p$
 the set $Q_{i_k\ldots i_{p}}$ is also non-connected and
 contains two components $Q^{(k)}_1$, $Q^{(k)}_2$ satisfying
 $\varphi_{i_1\ldots i_{k-1}}(Q_i^{(k)})
 \cap\partial \tilde{K}=Q^{(1)}_i$ also joined in
 $\partial \tilde{K}$ by an arc $\Delta^{(k)}$ so that
 $\varphi_{i_1\ldots i_{k-1}}(\Delta^{(k)})
 \cap\partial \tilde{K}=\Delta^{(1)}$.

  Let $N(\Delta^{(k)}$ be the set of all unit outer normal
  vectors to $\partial\tilde K$
 at points $x\in \Delta^{(k)}$. All these sets are open arcs. Since $\varphi_i$
  preserve angles between vectors, these arcs
  have  the same length.
   They cannot coincide for $k_1\neq k_2$ because
 the diameters of corresponding sets $\Delta^{(k_1)}$ and $\Delta^{(k_2)}$
 are different. They cannot have non-empty intersection, because of
 the OCSC condition. Therefore there cannot be more
 than $\frac{2\pi}{\theta}$ of such arcs.

 So, if $p_1>\frac{2\pi}{\theta}$, then each of the sets
 $Q_{i_1\ldots
i_{p_1}}$ is connected.

 Therefore the $p_1$-th refinement of the system $S$ satisfies
 conditions  1) and OCSC.

\vspace{10pt}

 A system $(\varphi_i,\ldots,\varphi_n)$ of contraction maps satisfying
conditions
 (\ref{e1}), (\ref{e2}), (\ref{e3}),and
  OCSC, will be called  {\it regular}.

 \begin{lem}
 If the system  $(\varphi_i,\ldots,\varphi_n)$ is regular,
 then the number of all components
 $Q_{i_1\ldots i_p}\subset\partial\tilde{K}$
 of order $p$, having non-empty interior in
 $\partial\tilde{K}$ does not exceed  $n^p$.
 \end{lem}

 {\it Proof.} The statement is obvious for $p=1$ .

 Suppose it's true for all components  $Q_{i_1\ldots i_{p-1}}$
 of order  $p-1$ .

 Suppose a component $Q_{i_1\ldots i_{p-1}}$ contains a
 $p$-component
 $Q_{i_1\ldots i_{p-1}j}$ different from $Q_{i_1\ldots i_{p-1}}$.
Then the endpoints   $\xi_j$, $\eta_j$ of the component $Q_j$
must satisfy either  $\varphi_{i_1\ldots i_{p-1}}(\xi_j) \in
 \stackrel{\circ}{Q}_{i_1\ldots i_{p-1}}$  or
$\varphi_{i_1\ldots i_{p-1}}
 (\eta_j) \in \stackrel{\circ}{Q}_{i_1\ldots i_{p-1}}$.

 Using this observation we can estimate the number of all  components of
  order
 p    $Q_{i_1\ldots i_p}$ having non-empty interior in
 $\partial\tilde{K}$.

 Consider the sets of unit normal vectors  $N^\circ_{i_1\ldots i_p}=
 \bigcup\limits_{x\in \stackrel{\circ}{Q}_{i_1\ldots i_p}}N_x$,
  where $\stackrel{\circ}{Q}_{i_1\ldots i_p}$ is the interior of
  the component  $Q_{i_1\ldots i_p}$ in $\partial\tilde{K}$.
If the set  $\stackrel{\circ}{Q}_{i_1\ldots i_p}$ is non-empty,
then $N^\circ_{i_1\ldots i_p}$ is an open arc of the unit circle.
It is clear that
 $N^\circ_{i_1}\supset
 N^\circ_{i_1i_2}\supset\ldots \supset N^\circ_{i_1\ldots i_p}\supset\ldots$
 and
 $N^\circ_{i_1\ldots i_p}=\varphi_{i_1\ldots i_p}(N^\circ_{i_p})\cap
 N^\circ_{i_1\ldots i_{p-1}}$.

If the set $N^\circ_{i_1\ldots i_p}$ differs from
 $N^\circ_{i_1\ldots i_{p-1}}$ one of the endpoints of the arc
 $N^\circ_{i_1\ldots i_p}$ must lie in $N^\circ_{i_1\ldots i_{p-1}}$.

 So, if $\beta^-_{i_p}$ and $\beta^+_{i_p}$ are the endpoints
 of the arc $N^\circ_{i_p}$, then one of the following
 inequalities hold:
 $$\beta^-_{i_1\ldots i_p}<\alpha_{i_1}+\ldots +\alpha_{i_{p-1}}+\beta^-_{i_p}<
 \beta^+_{i_1\ldots i_p}$$
 $$\beta^-_{i_1\ldots i_p}<\alpha_{i_1}+\ldots +\alpha_{i_{p-1}}+\beta^+_{i_p}<
 \beta^+_{i_1\ldots i_p}$$

 The sum $\alpha_{i_1}+\ldots
 +\alpha_{i_{p-1}}$ is the same for different
  permutations of $i_1\ldots i_{p-1}$, whereas the sets
$N^\circ_{i_1\ldots i_{p-1}}$
 are disjoint. Therefore we can take the union of all those
 $N^\circ_{i_1\ldots i_{p-1}}$ for which
 $\alpha_{i_1}+\ldots
 +\alpha_{i_{p-1}}=\gamma$, and denote it by $N^\circ_\gamma$.
 The union of the interiors  $\stackrel{\circ}{Q}_{i_1\ldots i_{p-1}}$
 of the corresponding $p-1$-components in $\partial\tilde{K}$ we denote by
 $Q^\circ_\gamma$. Let
 $p_\gamma$ be the number of these components.

 Thus, if a component $\stackrel{\circ}{Q}_{i_1\ldots i_p}\subset P_\gamma$
 is different from $\stackrel{\circ}{Q}_{i_1\ldots i_{p-1}}$,
 then one of the conditions  $\gamma+\beta^-_{i_p}\in
 N^\circ_\gamma$,
 $\gamma+\beta^+_{i_p}\in N^\circ_\gamma$ must hold.

 More exactly, if $p-1$-component
 $\stackrel{\circ}{Q}_{i_1\ldots i_{p-1}}$
contains $q+1$\ \   $p$-components
   $Q_{i_1\ldots i_{p}}$,...., $Q_{i_1\ldots
   i_{p-1}(i_{p}+q)}$,  then the system of inequalities
 $$\left\{\begin{array}{ccc}
   \beta^-_{i_1\ldots i_{p-1}}<&\alpha_{i_1}+\ldots +\alpha_{i_{p-1}}+
     \beta^+_{i_p+k}&<\beta^+_{i_1\ldots i_{p-1}} \\
   \beta^-_{i_1\ldots i_{p-1}}<&\alpha_{i_1}+\ldots +\alpha_{i_{p-1}}+
     \beta^-_{(i_p+k+1)}&<\beta^+_{i_1\ldots i_{p-1}}
   \end{array}\right.$$ holds for $k=0,...,q-1$.

Therefore if $m$ is a total number of those angles $\beta^+_j$ and
$\beta^-_j$, for which the conditions $\gamma+\beta^-_j\in
N^\circ_\gamma$ or
 $\gamma+\beta^+_j\in N^\circ_\gamma$ hold,
 then the number of components of order $p$ contained in
  $Q_\gamma$ is not greater than $p_\gamma+m$.

 The number of different sets $Q_\gamma$,
is in it's turn no greater than the number of different summands
in the expansion of
 $(x_1+\ldots +x_n)^{p-1}$, which is equal to
 $\frac{(n+p-2)!}{(p-1)!(n-1)!}$. Therefore the number of
different components ( having in $\partial\tilde{K}$ non-empty
interior) of  order $p$ exceeds the number of order $(p-1)$
components by a number not greater than
$n\cdot\frac{(n+p-2)!}{(p-1)!(n-1)!}$.

 So the total number of components of order $p$ is less or equal
 to
 $$n\left(1+\frac{n!}{1!(n-1)!}+\frac{(n+1)!}{2!(n-1)!}+\ldots+
 \frac{(n+p-2)!}{(p-1)!(n-1)!}\right)=\frac{(n+p-1)!}{(p-1)!(n-1)!},$$
 which is less than $p^{n+1}$.

 \begin{teo}
 If the system $(\varphi_1\ldots \varphi_n)$ satisfies the
OCS condition,  then Hausdorff dimension of the set $F$ of the
vertices of
 $\partial\tilde{K}$ is zero.
 \end{teo}

 {\it Proof.} We can suppose the system
 $(\varphi_1\ldots\varphi_n)$ is regular.
 The union of all $Q_{i_1\ldots i_p}$ consisting of one point is
 at most countable, so we may consider the set
  $F'=F\setminus\cup Q_{i_1\ldots i_p}$,
   which has the same Hausdorff dimension as  $F$.

 For each $p$,
 the set $F'$ is covered by components $Q_{i_1\ldots i_p}$,
  having non-empty interior in $\partial\tilde{K}$.
  The total number of such components does not exceed
 $p^{n+1}$, and  the diameter of each of them is less or equal to
  $q^p\cdot diam(\tilde{K})$. Since
 $\lim\limits_{p\rightarrow\infty}\left(-\frac{\ln p^{n+1}}{\ln q^p}\right)=0$,
 Hausdorff dimension of the set $F'$,
  and,therefore of $F$, is zero \cite{Fal}.

 \begin{sled}\label{sz}
If the system $\varphi_1,\ldots,\varphi_n$ satisfies
 OCSC and all the angles $\alpha_1,\ldots$ $,\alpha_n$ are
 commensurable with
 $\pi$, then the set  $\tilde{K}$ is a convex finite polygon.
 \end{sled}

 {\it Proof.}  Since the number of sides of $\tilde{K}$ ( Theorem \ref{t_1})
is finite, the set $F$ has finite number of components.
 All they have zero measure therefore each of them is a point.


\begin{thebibliography}{2222}

\bibitem{Fal}
{\em Falconer K.J.}: Fractal geometry:mathematical foundations and
applications . --  J.Wiley and Sons, New York, 1990.

\bibitem{Hut}
{\em Hutchinson J.}: Fractals and self-similarity. -- Indiana
Univ. Math. J., 30, No 5, 1981, pp.713-747.


\end{thebibliography}
 \end{document}